
\magnification=1100
\overfullrule0pt

\input amssym.def


\def\qed{\hbox{\hskip 1pt\vrule width4pt height 6pt depth1.5pt \hskip 1pt}}

\def\CC{{\Bbb C}}
\def\RR{{\Bbb R}}
\def\ZZ{{\Bbb Z}}

\def\cB{{\cal B}}

\def\cF{{\cal F}}

\def\supp{{\rm supp}}




\font\smallcaps=cmcsc10
\font\titlefont=cmr10 scaled \magstep1

\font\sectionfont=cmbx10
\font\tinyrm=cmr10 at 8pt


\newcount\sectno
\newcount\subsectno
\newcount\resultno

\def\section #1. #2\par{
\sectno=#1
\resultno=0
\bigskip\noindent{\sectionfont #1.  #2}~\medbreak}

\def\subsection #1\par{\bigskip\noindent{\it  #1} \medbreak}


\def\prop{ \global\advance\resultno by 1
\bigskip\noindent{\bf Proposition \the\sectno.\the\resultno. }\sl}
\def\lemma{ \global\advance\resultno by 1
\bigskip\noindent{\bf Lemma \the\sectno.\the\resultno. }
\sl}
\def\remark{ \global\advance\resultno by 1
\bigskip\noindent{\bf Remark \the\sectno.\the\resultno. }}
\def\example{ \global\advance\resultno by 1
\bigskip\noindent{\bf Example \the\sectno.\the\resultno. }}
\def\cor{ \global\advance\resultno by 1
\bigskip\noindent{\bf Corollary \the\sectno.\the\resultno. }\sl}
\def\thm{ \global\advance\resultno by 1
\bigskip\noindent{\bf Theorem \the\sectno.\the\resultno. }\sl}
\def\defn{ \global\advance\resultno by 1
\bigskip\noindent{\it Definition \the\sectno.\the\resultno. }\slrm}
\def\endthm{\rm\bigskip}

\def\endprop{\rm\bigskip}

\def\pf{\rm\bigskip\noindent{\it Proof. }}
\def\endpf{\qed\hfil\bigskip}


\def\formula{\global\advance\resultno by 1
\eqno{(\the\sectno.\the\resultno)}}
\def\formulano{\global\advance\resultno by 1 (\the\sectno.\the\resultno)}
\def\tableno{\global\advance\resultno by 1
\the\sectno.\the\resultno. }
\def\lformula{\global\advance\resultno by 1
\leqno(\the\sectno.\the\resultno)}

\def\monthname {\ifcase\month\or January\or February\or March\or April\or
May\or June\or
July\or August\or September\or October\or November\or December\fi}

\newcount\mins  \newcount\hours  \hours=\time \mins=\time
\def\now{\divide\hours by60 \multiply\hours by60 \advance\mins by-\hours
     \divide\hours by60         
     \ifnum\hours>12 \advance\hours by-12
       \number\hours:\ifnum\mins<10 0\fi\number\mins\ P.M.\else
       \number\hours:\ifnum\mins<10 0\fi\number\mins\ A.M.\fi}


\nopagenumbers
\def\runningtitle{\smallcaps Calibrated representations}
\headline={\ifnum\pageno>1\eoheadline\else\firstheadline\fi}
\def\names{\smallcaps arun ram}
\def\firstheadline{\noindent Preliminary Draft \hfill  \today}
\def\firstheadline{}
\def\eoheadline{\ifodd\pageno\oddheadline\else\evenheadline\fi}
\def\oddheadline{\tenrm\hfil\runningtitle\hfil\folio}
\def\evenheadline{\tenrm \folio\hfil{\names}\hfil}

\vphantom{$ $}  
\vskip.75truein

\centerline{\titlefont Calibrated representations of affine Hecke algebras}
\bigskip
\centerline{\rm Arun Ram${}^\ast$}
\centerline{Department of Mathematics}
\centerline{Princeton University}
\centerline{Princeton, NJ 08544}
\centerline{{\tt rama@math.princeton.edu}}
\centerline{Preprint: June 14, 1998}

\footnote{}{\tinyrm ${}^\ast$ Research supported in part by National
Science Foundation grant DMS-9622985, and a Postdoctoral Fellowship
at Mathematical Sciences Research Institute.}

\bigskip

\noindent{\bf Abstract.}
This paper introduces the notion of calibrated representations 
for affine Hecke algebras and classifies and constructs all finite
dimensional irreducible calibrated representations.  The main results are that (1)
irreducible calibrated representations are indexed by placed skew shapes, 
(2)  the dimension of an irreducible calibrated representation is the number of
standard Young tableaux corresponding to the placed skew shape and
(3)  each irreducible calibrated representation is constructed explicitly
by formulas which describe the action of each generator of the affine Hecke algebra
on a specific basis in the representation space.  This construction is a
generalization of A. Young's seminormal construction of the irreducible
representations of the symmetric group.  In this sense Young's
construction has been generalized to arbitrary Lie type.

\section 0. Introduction

The affine Hecke algebra was introduced by Iwahori and Matsumoto [IM]
as a tool for studying the representations of a $p$-adic Lie group.
In some sense, all irreducible principal series representations of
the $p$-adic group can be determined by classifying the representations
of the corresponding affine Hecke algebra.  Unfortunately, it is not
so easy to determine the irreducible representations of the affine Hecke
algebra.

Kazhdan and Lusztig [KL] (see also the important work of Ginzburg [CG])
gave a geometric classification of the irreducible representations of
the affine Hecke algebra. This classification is a $q$-analogue of Springer's
construction of the irreducible representations of the Weyl group on the cohomology
of unipotent varieties.  In the $q$-case, K-theory takes the place
of cohomology and the irreducible representations of the affine
Hecke algebra are constructed
as quotients of the K-theory of special subvarieties of
the flag variety.  Although the classification of Kazhdan and Lusztig
is an incredible tour-de-force it is difficult to obtain combinatorial
information from this geometric construction.  For example, it is difficult
to determine the dimensions of the irreducible modules.

In this paper I give a new construction of a large family
of irreducible modules of the affine Hecke algebra.  
The basis vectors are labeled by
generalized standard Young tableaux and the action of each generator on each
basis element is given explicitly.
This construction is a generalization of Young's
seminormal construction of the irreducible representations of the symmetric group.
In order to obtain this generalization I have had to
generalize the concept of standard Young tableaux to arbitrary Lie type.

The modules which I construct I
have termed ``calibrated'' modules.  Specifically,
a calibrated module is a module which has a basis of simultaneous
eigenvectors for all the elements of a large commutative algebra
inside the affine Hecke algebra.
This is analogous to the situation which occurs for representations
of complex semisimple Lie algebras where every finite dimensional
module is a direct sum of its weight spaces.
In contrast to the complex semisimple Lie algebra case, it is {\it never} 
true that all irreducible representations of the
affine Hecke algebra are calibrated.   

The irreducible really calibrated modules for the affine Hecke algebra are indexed
by placed skew shapes, where ``placed skew shape'' is a  generalization of the
usual skew shape from  combinatorial representation theory and symmetric function
theory. As is to be expected, these new generalized skew shapes and standard 
Young tableaux reduce to the classical objects in the Type A case.
This reduction is given in [Ra2].

\subsection{Remarks on the results in this paper}

(1) It is quite a surprise that the seminormal construction
of A. Young fits so nicely into general Lie type.  Up to now,
the general feeling has been that Young's results are very special to the
symmetric group and the Type A case.  The direct generalization
of Young's seminormal construction to arbitrary Lie type
which is obtained in this paper shows that this is not the case at all.
The results of [Ra4] indicate how 
Young's natural basis can also be generalized to arbitrary Lie type. 

The seminormal representations of A. Young have been previously generalized
to Iwahori-Hecke algebras of Type A by Hoefsmit [H] and Wenzl [Wz] independently,
to Iwahori-Hecke algebras of types B and D by Hoefsmit [H] and to the cyclotomic
Hecke algebras $H_{r,1,n}$ by Ariki and Koike [AK].  All of these generalizations
use classical standard Young tableaux and similar formulas for the action
of the generators of the Hecke algebra.  Using certain surjective
homomorphisms [A] from the affine Hecke algebras of type A to the algebras
$H_{r,1,n}$ one can easily show that these earlier constructions are type A special
cases of the general type construction given in this paper.

In a previous paper [Ra1] I gave a method for generalizing Young's theory of
seminormal representations to general Lie type.  I now believe that this earlier
idea was not the ``proper'' way to proceed.  The method here is much more natural
and yields a cleaner and more beautiful theory.  Young's classical formulas
for the seminormal representations of the symmetric group $S_n$ work in general Lie
type with no change at all!!  The only previously missing ingredient was a good
general type definition of standard Young tableaux.

(2)  In the classical theory of representations of the symmetric group $S_n$
the ``skew shape representations'' are particularly well behaved $S_n$-modules.  On
the other hand there never seemed to be any a priori raison d'etre for
skew shape  representations via which one could generalize this concept to
Weyl groups and Iwahori-Hecke algebras of other Lie types.
The results in this paper show that skew shape representations
do arise in a perfectly natural way.  They correspond
to {\it irreducible representations of affine
Hecke algebras}.  In the type $A$ cases one recovers the classical skew shape
representations by restricting to the Iwahori-Hecke algebra inside the affine Hecke
algebra.

(3) The two main techniques used in this paper are generalizations
of the techniques of Matsumoto [Ma] and Rodier [Ro].  In particular,
the $\tau$-operators and the calibration
graphs $\Gamma(t)$ introduced in Section 2 are generalizations 
of the intertwining operators of Matsumoto and of the graphs used by
Rodier, respectively.  I have changed the role of
the intertwining operators by having them be ``left'' operators instead
of ``right'' operators.  This means that they are no longer
intertwining but there are other benefits to 
using these operators in this fashion.

(4)  Heckman and Opdam [HO1-2] introduced a new ``harmonic analysis'' approach
to the representations of the affine Hecke algebra.  In their work they
also used the sets $Z(t)$ and $P(t)$ which we use in section 3.
These sets arise naturally in their work as the
zeros and poles of a certain Harish-Chandra $c$-function.
In this paper these sets describe the behaviour of the $\tau$-operators
mentioned in remark (3).  This means that there is a strong connection between
the $c$-function and these operators.  The approach of Heckman-Opdam
becomes difficult when one needs to compute the residues of the
$c$-function at certain singular points.  In some cases these difficulties
can be surmounted by using the methods of this paper.

(5)  The Kazhdan-Lusztig construction of the irreducible representations
of the affine Hecke algebras shows that the structure of these representations
is intimately connected with the geometry of certain subvarieties $\cB_{s,u}$ of the
flag variety.  It is my hope that the methods which I have used here for 
studying representations of affine Hecke algebras, which are mostly
combinatorial in nature, will be useful for studying the geometry of the
$\cB_{s,u}$ varieties used in the Kazhdan-Lusztig construction.  I am hoping
that the standard Young tableaux introduced here can be used as index sets for
the connected components of these varieties.

\subsection{Acknowledgements}

This paper is the first in a series [Ra2-4] [RR1-2] on
representations of affine Hecke algebras.  I have benefited from
conversations with many people. To choose only a few, there were discussions with
S. Fomin, F. Knop, L. Solomon, M. Vazirani and N. Wallach which played an
important role in my progress.  There were several times when
I tapped into J. Stembridge's fountain of useful knowledge about root systems. 
D.-N. Verma helped at a crucial juncture by suggesting that I look at the paper of
Steinberg.  G. Benkart was a very patient listener on many occasions.  H. Barcelo,
P. Deligne, T. Halverson, R. Macpherson and R. Simion all gave large amounts of time
to let me tell them my story and every one of these sessions was helpful to me in
solidifying my understanding.

I single out Jacqui Ramagge with special thanks for everything she
has done to help with this project: from the most mundane typing and picture
drawing to deep intense mathematical conversations which helped to sort out
many pieces of this theory.  Her immense contribution is evident in
that some of the papers in this series on representations of affine Hecke
algebras are joint papers.  

A portion of this research was done during a semester stay at Mathematical
Sciences Research Institute where I was supported by a Postdoctoral
Fellowship.  I thank MSRI and National Science Foundation for support of
my research.

\bigskip

\section 1. The affine Hecke algebra

\medskip

Let $R$ be a reduced irreducible root system in $\RR^n$, fix a set of positive
roots $R^+$ and let $\{\alpha_1,\ldots, \alpha_n\}$ be the corresponding
simple roots in $R$.  Let $W$ be the Weyl group corresponding to $R$.
Let $s_i$ denote the simple reflection in $W$ corresponding to the simple
root $\alpha_i$ and recall that $W$ can be presented by 
generators $s_1,s_2,\ldots, s_n$ and relations
$$\matrix{
s_i^2&=&1, &&\hbox{for $1\le i\le n$,} \cr
\underbrace{s_is_js_i\cdots}_{m_{ij} {\rm \ \ factors}}
&=& \underbrace{s_js_is_j\cdots}_{m_{ij} {\rm \ \ factors}}\;,
&\qquad &\hbox{for $i\ne j$,}  \cr
}
$$
where $m_{ij} = 
\langle \alpha_i,\alpha_j^\vee\rangle \langle\alpha_j,\alpha_i^\vee\rangle$.

Fix $q\in \CC^*$ such that $q$ is not a root of unity.  
The {\it Iwahori-Hecke algebra}
$H$ is the associative algebra over $\CC$ defined by 
generators $T_1,T_2,\ldots, T_n$ and relations
$$\matrix{
T_i^2&=& (q-q^{-1})T_i+1, &&\hbox{for $1\le i\le n$,} \cr
\underbrace{T_iT_jT_i\cdots}_{m_{ij} {\rm \ \ factors}}
&=& \underbrace{T_jT_iT_j\cdots}_{m_{ij} {\rm \ \ factors}}\;,
&\qquad &\hbox{for $i\ne j$,}  \cr
}
\formula$$
where $m_{ij}$ are the same as in the presentation of $W$.
For $w\in W$ define $T_w=T_{i_1}\cdots T_{i_p}$ where
$s_{i_1}\cdots s_{i_p} = w$ is a reduced
expression for $w$.  By [Bou, Ch. IV \S 2 Ex. 23], the element
$T_w$ does not depend on the choice of the reduced expression.
The algebra $H$ has dimension $|W|$ and the set $\{T_w\}_{w\in W}$
is a basis of $H$.

The {\it fundamental weights} are the elements $\omega_1,\ldots, \omega_n$
of $\RR^n$ given by
$$\langle \omega_i, \alpha_j^\vee\rangle = \delta_{ij},
\qquad\hbox{where}\quad
\alpha_i^\vee ={2\alpha_i\over \langle \alpha_i,\alpha_i\rangle}
$$
and $\delta_{ij}$ is the Kronecker delta.
The {\it weight lattice} is the $W$-invariant lattice in $\RR^n$ given by
$$P= \sum_{i=1}^n \ZZ\omega_i.$$
Let $X$ be the abelian group $P$ except written multiplicatively.  In
other words,
$$X = \{ X^\lambda \ |\ \lambda\in P\},
\quad\hbox{and}\quad
X^\lambda X^\mu = X^{\lambda+\mu} = X^\mu X^\lambda,
\quad\hbox{for $\lambda,\mu\in P$}.$$
Let $\CC[X]$ denote the group algebra of $X$.
There is a $W$-action on $X$ given by
$$wX^\lambda = X^{w\lambda}
\qquad\hbox{for $w\in W$, $X^\lambda\in X$,}$$
which we extend linearly to  a $W$-action on 
$\CC[X]$.

The {\it affine Hecke algebra} $\tilde H$ associated to $R$ and $P$ is the
algebra given by
$$\tilde H  
= \CC\hbox{-span} \{ T_w X^\lambda \ |\ w\in W, X^\lambda\in X\}$$ 
where the multiplication of the $T_w$ is as in the Iwahori-Hecke algebra
$H$, the multiplication of the $X^\lambda$ is as in $\CC[X]$ and
we impose the relation
$$X^\lambda T_i = T_i X^{s_i\lambda} + (q-q^{-1})
{X^\lambda-X^{s_i\lambda}\over 1-X^{-\alpha_i}}, 
\qquad\hbox{for $1\le i\le n$ and $X^\lambda\in X$.}
\formula$$
This formulation of the definition of $\tilde H$ is due to Lusztig [Lu2]
following work of Bernstein and Zelevinsky.
The elements $T_wX^\lambda$, $w\in W$, $X^\lambda\in X$, form a basis
of $\tilde H$.  

\thm (Bernstein, Zelevinsky, Lusztig [Lu1, 8.1])
The center of $\tilde H$ is  
$\CC[X]^W = \{ f\in \CC[X] \ |\ wf=f\}$.
\endthm


\section 2. Weight spaces and calibration graphs

\subsection{Weights}

\smallskip
Let 
$$T=\{ \hbox{group homomorphisms $t\colon X \to \CC^*$}\}.$$
The {\it torus} $T$ is an abelian group with a $W$-action given by
$(wt)(X^\lambda) = t(X^{w^{-1}\lambda})$.  Any element $t\in T$ is determined
by the values $t(X^{\omega_1}), t(X^{\omega_2}), \ldots, t(X^{\omega_n})$.
For any element $t\in T$ define the {\it polar decomposition} 
$$t=t_rt_c,\qquad\hbox{ $t_r,t_c\in T$ such that $t_r(X^\lambda)\in \RR_{>0}$,
and $|t_c(X^\lambda)|=1$,}$$
for all $X^\lambda\in X$.  There is a unique $\mu\in \RR^n$ and a unique 
$\nu\in \RR^n/P$ such that
$$t_r(X^\lambda)=e^{\langle \mu,\lambda\rangle}
\qquad\hbox{and}\qquad
t_c(X^\lambda)=e^{2\pi i\langle \nu,\lambda\rangle},
\qquad\hbox{for all $\lambda\in P$.}\formula$$
In this way we identify the sets $T_r=\{ t\in T\ |\ t=t_r\}$
and $T_c=\{ t\in T\ |\ t=t_c\}$ with $\RR^n$
and $\RR^n/P$, respectively.

\medskip
\subsection{Weight spaces}
\medskip

Let $M$ be a finite dimensional $\tilde H$-module.
For each $t\in T$ the {\it $t$-weight space of $M$} 
and the {\it generalized $t$-weight space} are the subspaces
$$\eqalign{
M_t &= \{ m\in M \ |\ X^\lambda m = t(X^\lambda)m
\hbox{\ for all $X^\lambda\in X$}\}
\qquad\hbox{and}  \cr
\cr
M_t^{\rm gen} &= \{ m\in M\ |\
\hbox{for each $X^\lambda\in X$,
$(X^\lambda-t(X^\lambda))^k m =0$ for some $k\in \ZZ_{>0}$}
\},  \cr
}
$$
respectively.
If $M_t^{\rm gen}\ne 0$ then $M_t\ne 0$.
In general $M\neq\bigoplus_{t\in T} M_t$, but we do have
$$
M=\bigoplus_{t\in T} M_t^{\rm gen}.
$$
This is a decomposition of $M$ into Jordan blocks
for the action of $\CC[X]$.  Define the {\it support} of $M$ to be
$$\supp(M) = \{ t\in T \ |\ M_t^{\rm gen}\ne 0\}.\formula$$  


\subsection{Principal series modules}

Let $t\in T$ and let $\CC v_t$ be the one dimensional
$\CC[X]$-module corresponding to the character 
$t\colon X\to \CC^*$. 
Specifically, $\CC v_t$ is the one dimensional vector space with
basis $\{v_t\}$ and $\CC[X]$-action given by 
$$X^\lambda v_t =t(X^\lambda)v_t,
\qquad\hbox{for all $X^\lambda\in X$.}
$$  
The {\it principal series representation} corresponding to $t$ is
$$
M(t) = \tilde H \otimes_{\CC[X]} \CC v_t.
\formula$$
The set $\{ T_w\otimes v_t \ |\ w\in W  \}$ is a basis for the $\tilde H$-module
$M(t)$ and $\dim(M(t))=|W|$.

If $w\in W$ and $X^\lambda\in X$ then the defining relation (1.2) for
$\tilde H$ implies that
$$X^\lambda (T_w\otimes v_t)= t(X^{w\lambda}) (T_w\otimes v_t)
+\sum_{u<w} a_u (T_u\otimes v_t),$$
where the sum is over $u<w$ in the Bruhat-Chevalley order and $a_u\in \CC$. 
It follows that the eigenvalues of $X$ on $M(t)$ are of the form $wt$, $w\in W$,
and by counting the multiplicity of each eigenvalue we have
$$M(t) = \bigoplus_{wt\in Wt} M(t)_{wt}^{\rm gen}
\qquad\hbox{where}\quad\dim(M(t)_{wt}^{\rm gen}) = |W_t|,
\quad\hbox{for all $w\in W$.}
\formula$$

\thm (Kato's irreducibility criterion [Ka])
Let $t\in T$ and define $P(t) = \{\alpha>0\ |\ t(X^\alpha) = q^{\pm 2} \}$.
The principal series module $M(t)$ is irreducible if and only if 
$P(t)=\emptyset$.
\endthm
 
\medskip
\noindent
{\bf Remark.}  Kato actually proves a more general result and thus needs
a further condition for irreducibility.  We have
simplified matters by specifying the weight lattice $P$ in our construction of the
affine Hecke algebra.  One can use any $W$-invariant lattice in
$\RR^n$ and Kato works in this more general situation.  When the
one uses the weight lattice $P$, a result of Steinberg [St, 4.2, 5.3] says that the
stabilizer
$W_t$ of a point $t\in T$ under the action of $W$ is always a reflection group.
Because of this Kato's criterion takes a simpler form.

\subsection{Irreducible modules}

\prop
Let $M$ be a finite dimensional $\tilde H$-module.
\smallskip
\item{(a)} For some $t\in T$, $M_t$ is nonzero.
\smallskip
\item{(b)} If $M$ is irreducible and $M_t\ne 0$ then $M$ is a
quotient of $M(t)$.
\smallskip
\item{(c)} If $M$ is irreducible then $\dim(M)\le |W|$.
\pf
(a)  As an $X(T)$-module
$M$ contains a simple submodule and this submodule must be
one-dimensional since all irreducible representations of a commutative
algebra are one-dimensional.  Thus, there is a nonzero weight vector in $M$.

\smallskip\noindent
(b)  Let $m_t$ be a nonzero vector in $M_t$.  Then there is a
unique $\tilde H$-module homomorphism determined by
$$
\matrix{
\phi \colon &M(t) &\longrightarrow &M \cr
&v_t &\longmapsto &m_t \cr
}
$$
where $v_t$ is as in the construction of $M(t)$ in (2.3)
This map is surjective since $M$ is irreducible.
Thus $M$ is a quotient of $M(t)$.

\smallskip\noindent
(c) follows from (b) since $\dim(M(t))=|W|$.
\endpf

\noindent It follows from Proposition (2.6b) and (2.4) that the support $\supp(M)$
of an irreducible $\tilde H$-module $M$ is contained in a single 
Weyl group orbit in $T$. Since $M$ is irreducible and $\tilde H$ has countable
dimension, Dixmier's version of Schur's lemma implies that $Z(\tilde H)$ acts on
$M$ by scalars.  Let $t\in T$ be such that
$$pM = t(p)M,
\qquad\hbox{for all $p\in Z(\tilde H)$.}
$$
Since $Z(\tilde H)=\CC[X(T)]^W$ it follows that
$t(p(X))=(wt)(p(X))$ for all $w\in W$.  The $W$-orbit $Wt$ of
$t$ is the {\it central character} of $M$.  We shall often abuse notation
and refer to any weight $s\in Wt$ as ``the central character'' of $M$.

\medskip\noindent
{\bf Remarks.}
\item{(a)}  The algebra $\CC[X]^W$ is the polynomial ring
$$\CC[X]^W= \CC[\chi^{\omega_1},\ldots,\chi^{\omega_n}],$$
where $\omega_1,\ldots,\omega_n$ are the fundamental weights in $P$
and $\chi^{\omega_1},\ldots, \chi^{\omega_n}$, are the corresponding 
Weyl characters. Thus, in order to specify the $W$-orbit of a weight
$t\in T$ 
it is sufficient to specify the $n$ complex numbers 
$t(\chi^{\omega_1}), t(\chi^{\omega_2}), \ldots, t(\chi^{\omega_n}).$ 
\item{(b)}  If $t\in T$ is real, i.e. $t=t_r$, then every element of $Wt$ is
also real.  Then there is a unique dominant $\gamma\in \RR^n$ such that the
element $s\in T$ given by
$$s(X^\lambda)=e^{\langle \gamma,\lambda\rangle},
\qquad\hbox{for all $\lambda\in P$,}$$
is in $Wt$ and this element is a canonical representative of the $W$-orbit $Wt$.  
In this way the dominant elements of $\RR^n$ index the real central characters.

\subsection{The $\tau$ operators}

The maps $\tau_i\colon M_t^{\rm gen}\to M_{s_it}^{\rm gen}$ 
defined below are local operators on $M$ in the sense that they
act on each weight space $M_t^{\rm gen}$ of $M$ separately.  The operator
$\tau_i$ is only defined on weight spaces $M_t^{\rm gen}$ such that
$t(X^{\alpha_i})\ne 1$.  

\prop
Let $t\in T$ such that $t(X^{\alpha_i})\ne 1$ and let $M$
be a finite dimensional $\tilde H$-module.  Define
$$
\matrix{
\tau_i \colon 
&M_t^{\rm gen} &\longrightarrow &M_{s_it}^{\rm gen} \cr
\cr
& m & \longmapsto &
\displaystyle{
\left(T_i - {q-q^{-1}\over 1-X^{-\alpha_i} }\right) m.} \cr
}
$$
\item{(a)}
The map $\tau_i\colon M_t^{\rm gen} \longrightarrow M_{s_it}^{\rm gen}$
is well defined.
\smallskip
\item{(b)}  As operators on $M_t^{\rm gen}$, \quad
$\displaystyle{ X^\lambda \tau_i = \tau_i X^{s_i\lambda}}$,
for all $X^\lambda\in X$.
\smallskip
\item{(c)}  As operators on $M_t^{\rm gen}$, \quad
$\displaystyle{
\tau_i\tau_i=
{(q-q^{-1}X^{\alpha_i})(q-q^{-1}X^{-\alpha_i})\over 
(1-X^{\alpha_i})(1-X^{-\alpha_i}) }.
}
$
\smallskip
\item{(d)}  Both maps $\tau_i\colon M_t^{\rm gen}\to M_{s_it}^{\rm gen}$
and $\tau_i\colon M_{s_it}^{\rm gen}\to M_t^{\rm gen}$
are invertible if and only if $t(X^{\alpha_i})\ne q^{\pm 2}$.
\smallskip
\item{(e)}  Let $1\le i\ne j\le n$ and let $m_{ij}$ be as in (1.1).
Then   
$$\underbrace{\tau_i\tau_j\tau_i\cdots}_{m_{ij} {\rm\ factors}}
= \underbrace{\tau_j\tau_i\tau_i\cdots}_{m_{ij} {\rm\ factors}},
$$
whenever both sides are well defined operators on $M_t^{\rm gen}$.
\pf
(a) 
Note that $(q-q^{-1})/(1-X^{-\alpha_i})$ is not a
well defined element of $\tilde H$ or $\CC[X]$ since it
is not a polynomial in $X^{-\alpha_i}$.  Because
of this we will be careful to view $(q-q^{-1})/(1-X^{-\alpha_i})$ 
only as a local operator.  Let us describe this operator more precisely.

The element $X^{\alpha_i}$ acts on $M_t^{\rm gen}$ by
$t(X^{\alpha_i})$ times a unipotent transformation.
As an operator on $M_t^{\rm gen}$, 
$1-X^{-\alpha_i}$ is invertible since it has determinant
$(1-t(X^{-\alpha_i}) )^d$ where $d=\dim(M_t^{\rm gen})$.
Since this determinant is nonzero $(q-q^{-1})/(1-X^{-\alpha_i})
=(q-q^{-1})(1-X^{-\alpha_i})^{-1}$ 
is a well defined operator on $M_t^{\rm gen}$.  Thus the definition of 
of $\tau_i$ makes sense.

The following calculation shows that
$\tau_i$ maps $M_t^{\rm gen}$ into $M_{s_it}^{\rm gen}$.
For the purposes of this calculation we are viewing all elements
of $\CC[X]$ as local operators on $M_t^{\rm gen}$ and
we abuse notation and denote the operator $(1-X^{-\alpha_i})^{-1}$
by $1/(1-X^{-\alpha_i})$.  We are able to do this without any problems
because $\CC[X]$ is commutative.
$$
\eqalign{
X^\lambda \tau_i m
&=\left(
T_iX^{s_i\lambda} + (q-q^{-1})
{X^\lambda-X^{s_i\lambda}\over 1-X^{-\alpha_i}}
-X^\lambda{q-q^{-1}\over 1-X^{-\alpha_i} }\right)m \cr
&=
\left(
T_iX^{s_i\lambda} - (q-q^{-1})
{X^{s_i\lambda}\over 1-X^{-\alpha_i} } \right)m \cr
&=
\left(
T_i - {q-q^{-1}\over 1-X^{-\alpha_i} } \right)
X^{s_i\lambda}m \cr
&= \tau_i X^{s_i\lambda}m.\cr
}
$$
(b) follows from the previous calculation.
\smallskip\noindent
(c)
If $t(X^{-\alpha_i})\ne 1$ then both 
$\tau_i\colon M_t^{\rm gen}\to M_{s_it}^{\rm gen}$ and
$\tau_i\colon M_{s_it}^{\rm gen}\to M_t^{\rm gen}$ are well 
defined.  If $m\in M_t^{\rm gen}$ then
$$\eqalign{
\tau_i\tau_i m
&=
\left(T_i - {q-q^{-1}\over 1-X^{-\alpha_i} }\right)
\left(T_i - {q-q^{-1}\over 1-X^{-\alpha_i} }\right) m \cr
&=
\left(T_i^2 - {q-q^{-1}\over 1-X^{-\alpha_i} }T_i
-T_i {q-q^{-1}\over 1-X^{-\alpha_i} }
+{(q-q^{-1})^2\over (1-X^{-\alpha_i})^2 }
\right) m \cr
&=
\left((q-q^{-1})T_i +1 - {q-q^{-1}\over 1-X^{-\alpha_i} }T_i
-{ q-q^{-1}\over 1-X^{\alpha_i} }T_i \right. \cr
&\qquad\qquad
\left.
-(q-q^{-1})^2
{ (1-X^{-\alpha_i})^{-1}-(1-X^{\alpha_i})^{-1}\over 1-X^{-\alpha_i} }
+{(q-q^{-1})^2\over (1-X^{-\alpha_i})^2 }
\right) m \cr
&=
\left((q-q^{-1})T_i +1 - (q-q^{-1})T_i
+(q-q^{-1})^2{-1-X^{-\alpha_i}+1\over (1-X^{-\alpha_i})^2}
\right) m \cr
&=
\left(1+{(q-q^{-1})^2\over (1-X^{-\alpha_i})(1-X^{\alpha_i})}
\right)m \cr
&=
 {2-X^{\alpha_i}-X^{-\alpha_i}+(q^2-2+q^{-2})\over 
(1-X^{\alpha_i})(1-X^{-\alpha_i}) }
 m \cr
&=
{ (q-q^{-1}X^{\alpha_i})(q-q^{-1}X^{-\alpha_i})\over 
(1-X^{\alpha_i})(1-X^{-\alpha_i}) } 
 m. \cr
}
$$
(d)  The operator $X^{\alpha_i}$ acts on 
$M_t^{\rm gen}$ as $t(X^{\alpha_i})$ times a unipotent transformation.
Similarly for $X^{-\alpha_i}$.  Thus, as an operator on $M_t^{\rm gen}$
$\det((q-q^{-1}X^{\alpha_i})(q-q^{-1}X^{-\alpha_i}))=0$
if and only if $t(X^{\alpha_i})=q^{\pm 2}$.
Thus part (c) implies that $\tau_i\tau_i$
is invertible if and only if $t(X^{\alpha_i})\ne q^{\pm 2}$.
The statement follows.
\smallskip\noindent
(e) Let us begin with a slight diversion which will be helpful in the proof. 
Let $t\in T$ be a generic element of $T$ and let $M(t)$ be the corresponding
principal series module. Since $t$ is generic,
$W_t=\{1\}$ and 
$$M(t)=\bigoplus_{w\in W} M(t)_{wt},
\qquad\hbox{and}\qquad \dim(M(t)_{wt})=1,$$
for all $w\in W$.  We have $M(t)_{wt}^{\rm gen}=M(t)_{wt}$ since $M(t)_{wt}$
is nonzero whenever $M(t)_{wt}^{\rm gen}$ is nonzero and 
we know that $\dim(M(t)_{wt}^{\rm gen})=1$.   
Let $w\in W$ such that $\ell(s_iw)=\ell(w)+1$.  
Since $t$ is generic $(wt)(X^{\alpha_i})\ne q^{\pm 2}$ for all $w\in W$. Thus by
part (d), the map 
$\tau_i\colon M(t)_{wt}^{\rm gen}\to M(t)_{s_iwt}^{\rm gen}$ 
is a bijection.  Using the vector $v_t\in M(t)_t$ and the maps $\tau_i$ we
can construct a basis $\{v_{wt}\}_{w\in W}$ of $M(t)$ given by
$$
v_{s_iwt}=\tau_iv_{wt},\qquad
\hbox{if $\ell(s_iw)=\ell(w)+1$.}
$$
This basis is uniquely determined by the conditions
$$\matrix{
\hbox{(\global\advance\resultno by 1
\the\sectno.\the\resultno a)} &X^\lambda v_{wt} = (wt)(X^\lambda)
v_{wt},\hfill &\quad &\hbox{for all $w\in W$ and $X^\lambda\in X(T)$,}\hfill
\cr
\cr
\hbox{(\the\sectno.\the\resultno b)}&\phantom{X^\lambda }v_{wt} =
\displaystyle{ T_w\otimes v_t + \sum_{u<w} a_{wu}(t)(T_u\otimes v_t),}
\hfil &&\hbox{where $a_{wu}(t)\in \CC$.}\hfill \cr
}$$

Now we proceed to the proof of the statement.
We may assume that $\tilde H$ is the affine Hecke algebra corresponding
to a rank two root system $R$ generated by simple roots $\alpha_i$ and 
$\alpha_j$.  Let $w_0$ be the longest element of
$W$.   Every element $w\in W$, $w\ne w_0$ has unique minimal length expression as a
product of generators
$s_i$ and $s_j$.  Let $T_w$ be the corresponding product of the $T_i$'s and
$T_j$'s.  Using the defining relation (1.2) for $\tilde H$ we expand to derive
$$\underbrace{\cdots \left(T_i-{q-q^{-1}\over 1-X^{\alpha_i} }\right)
\left(T_j-{q-q^{-1}\over 1-X^{\alpha_j} }\right)
\left(T_i-{q-q^{-1}\over 1-X^{\alpha_i} }\right)}_{m_{ij} {\rm\ factors}}
=\underbrace{\cdots T_iT_jT_i}_{m_{ij} {\rm\ factors}}+\sum_{w<w_0} T_{w}
P_{w},
\formula$$ 
where the sum is over $w\in W$ such that $w\ne w_0$ and
$P_w$ are rational functions of the $X^\alpha$, $\alpha\in R$.  Similarly,
$$\underbrace{\cdots \left(T_j-{q-q^{-1}\over 1-X^{\alpha_j} }\right)
\left(T_i-{q-q^{-1}\over 1-X^{\alpha_i} }\right)
\left(T_j-{q-q^{-1}\over 1-X^{\alpha_j} }\right)}_{m_{ij} {\rm\ factors}}
=\underbrace{\cdots T_jT_iT_j}_{m_{ij} {\rm\ factors}}+\sum_{w<w_0} T_w Q_w,
\formula$$
where, as before, the sum is over $w\in W$ such that $w\ne w_0$ and
$Q_w$ are rational functions of the $X^\alpha$, $\alpha\in R$.
We shall show that $P_w=Q_w$.  

Let $t\in T$ be generic and let $M(t)$ be the corresponding principal
series module for $\tilde H$.  By the analysis in the previous paragraph
we have
$$\eqalign{
v_{w_0t} = \underbrace{\cdots \tau_i\tau_j\tau_i}_{m_{ij} {\rm\ factors}}v_{t} 
&= \underbrace{\cdots T_iT_jT_i}_{m_{ij} {\rm\ factors}} v_{t} +
\sum_{w<w_0} T_w P_w v_{t} \cr
&=T_{w_0} \otimes v_{t} +
\sum_{w<w_0} t(P_w) T_w\otimes v_{t}.  \cr
}$$
and it follows from (2.8b) that $t(P_w)=a_{w_0 w}(t)$ for all $w\in W$, $w\ne w_0$.
One shows similarly that $t(Q_w)=a_{w_0 w}(t)$ for all $w\in W$, $w\ne w_0$.

We have shown that, for each $w\in W$, $t(P_w)=t(Q_w)$ for all generic $t\in T$.  
Since $P_w$ and $Q_w$ are rational functions which coincide on all generic
points it follows that 
$$P_w=Q_w
\qquad\hbox{for all $w\in W$, $w\ne w_0$.}
\formula$$
Thus,
$$\eqalign{
\underbrace{\cdots \tau_i\tau_j\tau_i}_{m_{ij} {\rm\ factors}}
&=\underbrace{\cdots \left(T_i-{q-q^{-1}\over 1-X^{\alpha_i} }\right)
\left(T_j-{q-q^{-1}\over 1-X^{\alpha_j} }\right)
\left(T_i-{q-q^{-1}\over 1-X^{\alpha_i} }\right)}_{m_{ij} {\rm\ factors}} \cr
&=\underbrace{\cdots \left(T_j-{q-q^{-1}\over 1-X^{\alpha_j} }\right)
\left(T_i-{q-q^{-1}\over 1-X^{\alpha_i} }\right)
\left(T_j-{q-q^{-1}\over 1-X^{\alpha_j} }\right)}_{m_{ij} {\rm\ factors}} 
= \underbrace{\cdots \tau_j\tau_i\tau_j}_{m_{ij} {\rm\ factors}}, \cr
}$$
whenever both sides are well defined operators on $M_t^{\rm gen}$.
\endpf

\vfill\eject

\subsection{The calibration graph}

Let $t\in T$.  Define a graph $\Gamma(t)$ with
$$
\matrix{
\hbox{Vertices:} \quad  &Wt,\hfil & \cr
\hbox{Edges:}  &wt \longleftrightarrow s_iwt,\quad
&\hbox{ if }\quad (wt)(X^{\alpha_i}) \ne q^{\pm2}. \cr
}
$$

\prop  If $M$ is a finite dimensional $\tilde H$-module then
$$\dim(M_t^{\rm gen})=\dim(M_{t'}^{\rm gen})$$
if $t$ and $t'$ are in the same connected component of the
calibration graph.
\pf
It follows from Proposition (2.7d) that if
there is an edge $wt\leftrightarrow s_iwt$ in
$\Gamma(t)$ then the map $\tau_i\colon M_{wt}^{\rm gen}\to M_{s_iwt}^{\rm
gen}$ is a bijection.  Thus, $\dim(M_{wt}^{\rm gen})=\dim(M_{s_iwt}^{\rm gen})$
if $t$ and $s_it$ are connected in the calibration graph $\Gamma(t)$.
\endpf

\cor If $M$ is an irreducible $\tilde H$-module with central
character $t$ then the support $\supp(M)$
is a union of connected components of the calibration graph $\Gamma(t)$.

\medskip
\subsection{The connected components of $\Gamma(t)$}

\medskip
Let $t\in T$ and define
$$
Z(t) = \{ \alpha>0\ |\  t(X^\alpha)=1 \}, \qquad\hbox{and}\qquad
P(t) =  \{ \alpha>0 \ |\ t(X^\alpha)=q^{\pm 2} \}.
$$
If $J\subseteq P(t)$ define
$$\cF^{(t,J)} = \{ w\in W \ |\ R(w)\cap Z(t) = \emptyset,\ \  R(w)\cap P(t) = J\},$$
where $R(w) = \{ \alpha>0 \ |\ w\alpha<0 \}$ is the {\it inversion set} of
$w$.
Define a {\it placed shape} to be a pair $(t,J)$ such that $t\in T$,
$J\subseteq P(t)$ and $\cF^{(t,J)}\ne \emptyset$.  The elements of the set
$\cF^{(t,J)}$ are called {\it standard tableaux} of shape $(t,J)$.

\thm  The connected components of the calibration graph $\Gamma(t)$ are given by 
the partition of the vertices according to the sets
$$\cF^{(t,J)}t,
\qquad\hbox{such that $J\subseteq P(t)$ and $\cF^{(t,J)}\ne \emptyset$.}
$$
\pf  Let us begin by introducing appropriate notation.
The chamber 
$$C=\{ x\in \RR^n \ |\ 
\hbox{$\langle x,\alpha\rangle >0$ for all $\alpha\in R^+$} \}$$
is a fundamental chamber for the action of $W$ on $\RR^n$ and the
complement $\RR^n\setminus (\bigcup_\alpha H_\alpha)$
of the hyperplanes
$$H_\alpha=\{ x\in \RR^n \ |\ \langle x,\alpha\rangle =0 \},
\qquad \alpha\in R^+,$$
in $\RR^n$ is the disjoint union of the chambers $\{w^{-1}C\ |\ w\in W\}$.  A
chamber $w^{-1}C$ is on the {\it positive side} of the hyperplane $H_\alpha$
if $\langle x,\alpha\rangle >0$ for all $x\in w^{-1}C$.
The chambers adjacent to $w^{-1}C$ are the
chambers $w^{-1}s_iC$, $1\le i\le n$, and the common face of $w^{-1}C$ and
$w^{-1}s_iC$ is contained in the hyperplane $H_{w^{-1}\alpha_i}$.  

Now let $t$ be as in the statement of the Theorem.  
A result of Steinberg [St, 3.15, 4.2, 5.3] says that
the stabilizer of $t$ is 
$$W_t = \langle s_\alpha \ |\ \alpha\in Z(t)\rangle,$$ 
the subgroup of $W$ generated by the reflections in the
hyperplanes orthogonal to the roots in $Z(t)$.  The elements
of the orbit $Wt$ can be identified with the cosets in 
$W/W_t$ and these can be identified with the chambers of 
$\RR^n\setminus (\bigcup_\alpha H_\alpha)$ which are on the positive
side of the hyperplanes $H_\alpha$, $\alpha\in Z(t)$.  Under this
bijection the element $wt\in Wt$ is identified with the chamber
$w^{-1}C$.  The elements $wt$ and $s_iwt$ are not connected by an edge 
in $\Gamma(t)$ if and only if the hyperplane
$H_{w^{-1}\alpha_i}$ containing the common face of the corresponding
(adjacent) chambers $w^{-1}C$ and $w^{-1}s_iC$ is the hyperplane 
$H_\beta$ for some root $\beta\in P(t)$.
In this way we can identify the graph $\Gamma(t)$ with
the graph with
$$
\matrix{
\hbox{Vertices:} \quad  &\hbox{chambers of $\RR^n\setminus (\bigcup_\alpha
H_\alpha)$ which are on the positive side of the}\hfill &\cr
&\hbox{\qquad  hyperplanes $H_\alpha$, $\alpha\in Z(t)$} ,\hfill & \cr
\hbox{Edges:}  &\hbox{faces of the chambers
which are not contained in the}\hfill &\cr
&\hbox{\qquad hyperplanes $H_\beta$, $\beta\in P(t)$.}
\hfill
\cr
}
$$
The (closures of the) chambers $w^{-1}C$ which are on the positive side of the
hyperplanes $H_\alpha$, $\alpha\in Z(t)$, form a convex region in $\RR^n$.
This region is a disjoint union of smaller convex regions bounded by the
hyperplanes $H_\beta$, $\beta\in P(t)$.  Each of the smaller
regions is on the positive side of some of the hyperplanes
$H_\beta$, $\beta\in P(t)$, and on the negative side of others.
In fact, it is determined by the set
$J\subseteq P(t)$ such that it is on
negative side of the hyperplanes $H_\beta$, $\beta\in J$. 
These smaller convex regions correspond to the
connected components of $\Gamma(t)$ and thus the connected components of
$\Gamma(t)$ are given by the sets $\{wt\ |\ w\in \cF^{(t,J)}\}$ where $J\subseteq
P(t)$ and 
$$\cF^{(t,J)} = \left\{ w\in  W \ \Bigg|\ 
\eqalign{ \hbox{$w^{-1}C$ is}\cr\cr\cr}
\eqalign{
&\hbox{on the positive
side of the hyperplanes $H_\alpha$, $\alpha\in Z(t)$,
} \cr
&\hbox{on the positive
side of the hyperplanes $H_\alpha$, $\alpha\in P(t)\setminus J$,
} \cr
&\hbox{on the negative side of the hyperplanes $H_\beta$, $\beta\in J$} \cr}
\right\}.$$
Since the chamber $w^{-1}C$ is on the positive side of a hyperplane
$H_\alpha$ if and only if $w\alpha > 0$ it follows that 
$$\cF^{(t,J)} = \{ w\in W \ |\ R(w)\cap Z(t) = \emptyset,\ \ R(w)\cap P(t) = J\}.
\qquad\hbox{\qed}$$


\bigskip
\section 3.  An $\tilde H$-module construction

\medskip

Let $\alpha_i$ and $\alpha_j$ be simple roots in $R$ and let
$R_{ij}$ be the rank two root subsystem of $R$ which is generated
by $\alpha_i$ and $\alpha_j$.  Let $W_{ij}$ be the Weyl group of
$R_{ij}$, the subgroup of $W$ generated by the simple reflections
$s_i$ and $s_j$.  A weight $t\in T$ is {\it calibratable} for $R_{ij}$ if
one of the following two conditions holds:
\smallskip
\itemitem{(a)} $t(X^{\alpha})\ne 1$ for all $\alpha\in R_{ij}$,
\smallskip
\itemitem{(b)} $R_{ij}$ is of type $C_2$ or $G_2$ (we may assume that
$\alpha_i$ is the long root and $\alpha_j$ is the short root),
$ut(X^{\alpha_i})=q^2$ and $ut(X^{\alpha_j})=1$ for some $u\in W_{ij}$,
and $t(X^{\alpha_i})\ne 1$ and  $t(X^{\alpha_j})\ne 1$.
\smallskip\noindent
A {\it placed skew shape} is a placed shape $(t,J)$ such that 
for all $w\in \cF^{(t,J)}$ and all pairs $\alpha_i, \alpha_j$ of
simple roots in $R$ the weight $wt$ is calibratable for $R_{ij}$.

\thm  Let $(t,J)$ be a placed skew shape and let $\cF^{(t,J)}$ be the
set of standard tableaux of shape $(t,J)$.  Define 
$$\tilde H^{(t,J)} =
\CC\hbox{-span} \{ v_w \ |\ w\in {\cF}^{(t,J)} \},$$
so that the symbols $v_w$ are a labeled basis of the
vector space $\tilde H^{(t,J)}$.
Then the following formulas make $\tilde H^{(t,J)}$
into an irreducible $\tilde H$-module:  For each $w\in {\cF}^{(t,J)}$,
$$
\matrix{
\hfill X^\lambda v_w &=& (wt)(X^\lambda) v_w, \hfill
&&\hbox{for $X^\lambda\in X$, and}\hfil \cr
\cr
\hfill T_i v_w 
&=& (T_i)_{ww} v_w + (q^{-1}+(T_i)_{ww}) v_{s_iw}, \hfill
&\qquad &\hbox{for $1\le i\le n$,}\hfil \cr
}
$$
where $\displaystyle{
(T_i)_{ww} = {q-q^{-1}\over 1 - (wt)(X^{-\alpha_i})}\;, }
$
and we set $v_{s_iw} = 0$ if $s_iw\not\in {\cF}^{(t,J)}$.
\endthm
\pf
Since $(t,J)$ is a placed skew shape $(wt)(X^{-\alpha_i})\ne 1$ for all $w\in
\cF^{(t,J)}$ and all simple roots $\alpha_i$.  This implies that the coefficient
$(T_i)_{ww}$ is well defined for all $i$ and $w\in \cF^{(t,J)}$.

By construction, the nonzero weight spaces of $\tilde
H^{(t,J)}$ are 
$(\tilde H^{(t,J)})_{wt}^{\rm gen}=(\tilde H^{(t,J)})_{wt}$ 
where $w\in \cF^{(t,J)}$.  These weight spaces have
dimension $1$ and are all in a single connected component of the calibration graph
$\Gamma(t)$.  If $N$ is a proper submodule of $\tilde H^{(t,J)}$ then we would
have $N_{wt}\ne 0$ and $N_{w't}=0$ for some $w\ne w'$, $w,w'\in \cF^{(t,J)}$.
But since $wt$ and $w't$ are in the same connected component of $\Gamma(t)$ this
would contradict Proposition 2.12.  Thus
$\tilde H^{(t,J)}$ is irreducible if it is an $\tilde H$-module.

It remains to show that the defining relations for $\tilde H$ are
satisfied.
\smallskip\noindent
(a) Let $w\in \cF^{(t,J)}$.  Then 
$$\eqalign{
\Big( X^{s_i\lambda}T_i+&(q-q^{-1})
{X^\lambda - X^{s_i\lambda}\over 1-X^{-\alpha_i} }
\Big)v_w \cr
&=\left(
(wt)(X^{s_i\lambda}){(q-q^{-1})\over 1-(wt)(X^{-\alpha_i})}
+
(q-q^{-1})
{(wt)(X^\lambda) - (wt)(X^{s_i\lambda})\over 1-(wt)(X^{-\alpha_i}) }
\right)v_w  \cr
&\qquad\qquad\quad
+ (s_iwt)(X^{s_i\lambda})(T_i)_{s_iw,w} v_{s_iw} \cr
&= {(q-q^{-1})\over 1-(wt)(X^{-\alpha_i}) } (wt)(X^\lambda) v_w
+ (T_i)_{s_iw,w} (wt)(X^\lambda) v_{s_iw} \cr
&= T_iX^\lambda v_w. \cr
}
$$
(b) Let $w\in \cF^{(t,J)}$.  Using the fact that
$(T_i)_{ww}+(T_i)_{s_iw,s_iw}=q-q^{-1}$ we have 
$$
\eqalign{
T_i^2v_w 
&= ((T_i)_{ww}^2 
+ (q^{-1}+(T_i)_{ww})(q^{-1}+(T_i)_{s_iw,s_iw}))v_w \cr
&\qquad\qquad\quad
+ (q^{-1}+(T_i)_{ww})((T_i)_{ww}+(T_i)_{s_iw,s_iw}) v_{s_iw} 
\cr
&= (T_i)_{ww}( (T_i)_{ww}+(T_i)_{s_iw,s_iw} ) v_w
+q^{-1}(q^{-1} + (T_i)_{ww}+(T_i)_{s_iw,s_iw}) )v_w \cr
&\qquad\qquad\quad
+ (q^{-1}+(T_i)_{ww}) (q-q^{-1}) v_{s_iw}\cr
&= (T_i)_{ww}(q-q^{-1}) v_w +(q^{-1}+(T_i)_{ww})(q-q^{-1})v_{s_iw}
+q^{-1}(q^{-1}+q-q^{-1}) v_w \cr
&= ((q-q^{-1})T_i +1)v_w.\cr
}
$$
(c)  The braid relation.  Let $\alpha_i$ and $\alpha_j$ be simple
roots in $R$ and let $w\in \cF^{(t,J)}$.  Since $(t,J)$ is a placed skew 
shape $wt$ is calibratable.  There are two distinct cases to consider.

\smallskip\noindent
{\it Case 1}:  When $wt$ is $R_{ij}$-regular, i.e.
$(uwt)(X^{\alpha_i})\ne 1$ and $(uwt)(X^{\alpha_j})\ne 1$ for all $u\in W_{ij}$.

Let us extend our notation $v_\sigma$ to all $\sigma\in W$ by assuming that
$v_\sigma=0$ whenever $\sigma\not\in \cF^{(t,J)}$.  Then, for all
$u\in W_{ij}$, the definition of the action of $T_i$ allows us to write
$$\eqalign{
\left(T_i-{(q-q^{-1})\over 1-(wt)(X^{-\alpha_i})}\right)v_{uw}
&=\left(q^{-1}+{(q-q^{-1})\over 1-(wt)(X^{-\alpha_i})}\right)
v_{s_iuw} \cr
&=\left({q-q^{-1}(wt)(X^{-\alpha_i})
\over 1-(wt)(X^{-\alpha_i})}\right) v_{s_iuw} \cr}
$$
whenever $\ell(s_iu)>\ell(u)$ and $(uwt)(X^{-\alpha_i})\ne 1$.
(This is correct because if $\ell(s_iu)>\ell(u)$ and $v_{uw}=0$ then
$v_{s_iuw}=0$.) 

Since
$wt$ is $R_{ij}$-regular all of the factors
in the following product are well defined and, by 
[Bou, Ch. VI \S 1 Cor. 2 to Prop. 17],
$$\eqalign{
&\prod_{\alpha\in R_{ij}^+} {q-q^{-1}(wt)(X^{-\alpha})
\over 1-(wt)(X^{-\alpha})} \cr
&\quad=\underbrace{\cdots 
\left({q-q^{-1}(s_js_iwt)(X^{-\alpha_i})\over
1-(s_js_iwt)(X^{-\alpha_i})}\right)
\left({q-q^{-1}(s_iwt)(X^{-\alpha_i})\over
1-(s_iwt)(X^{-\alpha_i})}\right)
\left({q-q^{-1}(wt)(X^{-\alpha_i})\over
1-(wt)(X^{-\alpha_i})}\right) }_{m_{ij} {\rm\ factors}} \cr
}
$$
where the product is over all positive roots in the root subsystem $R_{ij}$
spanned by $\alpha_i$ and $\alpha_j$.  Thus we get that
$$\eqalign{
&\left(\prod_{\alpha\in R_{ij}^+} {q-q^{-1}(wt)(X^{-\alpha})
\over 1-(wt)(X^{-\alpha})}\right)v_{w_0w} \cr
&
=\underbrace{\cdots 
\left(T_i-{q-q^{-1}\over 1-(s_js_iwt)(X^{-\alpha_i})}\right)
\left(T_j-{q-q^{-1}\over 1-(s_iwt)(X^{-\alpha_i})}\right)
\left(T_i-{q-q^{-1}\over 1-(wt)(X^{-\alpha_i})}\right)}_{m_{ij} {\rm\ factors}}
v_w
\cr
&
=\underbrace{\cdots \left(T_i-{q-q^{-1}\over 1-X^{\alpha_i}
}\right)
\left(T_j-{q-q^{-1}\over 1-X^{\alpha_j} }\right)
\left(T_i-{q-q^{-1}\over 1-X^{\alpha_i} }\right)}_{m_{ij} {\rm\ factors}}v_w \cr
&
=\underbrace{\cdots T_iT_jT_i}_{m_{ij} {\rm\ factors}} v_w +
\sum_{u<u_0} T_u P_u v_w,
\cr}
$$
where the notation in the last line is exactly the same as in (2.9).
By similar reasoning we obtain
$$\left(\prod_{\alpha\in R_{ij}^+} {q-q^{-1}(wt)(X^{-\alpha})
\over 1-(wt)(X^{-\alpha})}\right)v_{w_0w}
=\underbrace{\cdots T_jT_iT_j}_{m_{ij} {\rm\ factors}} v_w +
\sum_{u<u_0} T_u Q_u v_w,$$
where the $Q_u$ are as in (2.10).  By (2.11), $P_u=Q_u$ for all $u\ne u_0$
in $W_{ij}$ and thus
$$\underbrace{\cdots T_iT_jT_i}_{m_{ij} {\rm\ factors}} v_w
=\underbrace{\cdots T_jT_iT_j}_{m_{ij} {\rm\ factors}} v_w.$$

\smallskip\noindent
{\it Case 2}:  Let $u\in W_{ij}$ be of minimal length such that
$(uwt)(X^{\alpha_i})=q^2$ and $(uwt)(X^{\alpha_j})=1$.   The only possibilities
are the following.
\medskip\noindent
{\bf Type $C_2$:}
\smallskip\noindent
(1)  $u=s_i$:  Then
$$\matrix{
X^{\alpha_i}v_w=q^{-2}v_w, 
&\qquad &X^{\alpha_j}v_w=q^2v_w, \cr
T_iv_w = -q^{-1}v_w, 
&\qquad &T_jv_w=qv_w.\cr
}
$$
(2)  $u=s_is_j$:  Then
$$\matrix{
X^{\alpha_i}v_w=q^2v_w, 
&\qquad &X^{\alpha_j}v_w=q^{-2}v_w, \cr
T_iv_w = qv_w, 
&\qquad &T_jv_w=-q^{-1}v_w.\cr
}
$$
\medskip\noindent
{\bf Type $G_2$:}
\smallskip\noindent
(1)  $u=s_i$:  Then
$$\matrix{
X^{\alpha_i}v_w=q^{-2}v_w, 
&\qquad &X^{\alpha_2}v_w=q^2v_w, \cr
T_iv_w = -q^{-1}v_w, 
&\qquad &T_jv_w=qv_w.\cr
}
$$
(2)  $u=s_is_j$ or $u=s_is_js_i$.  Then both $w$ and $s_iw$ are
in $\cF^{(t,J)}$ and the action of
$X^{\alpha_i}$,  $X^{\alpha_j}$,
$T_i$ and $T_j$ on $\CC$-span$\{v_w,v_{s_iw}\}$
is given by the matrices:
$$\matrix{
X^{\alpha_i}=\pmatrix{ q^4 &0\cr 0 &q^{-2} \cr},
&\qquad &X^{\alpha_j}=\pmatrix{ q^{-4} &0 \cr 0 &q^2 \cr}, \cr
\cr
T_i= \pmatrix{
\displaystyle{ {q-q^{-1}\over 1-q^{-4}} } 
& \displaystyle{ {q-q^3\over 1-q^4} } \cr
\displaystyle{{q-q^{-5}\over 1-q^{-4}} } 
& \displaystyle{{q-q^{-1}\over 1-q^4} } \cr
},
&&T_j= \pmatrix{ 
-q^{-1}  &0  \cr
 0       &q         \cr
}.
\cr
}$$
(3)  $u=s_is_js_is_j$:  Then
$$\matrix{
X^{\alpha_i}v_w=q^2v_w, 
&\qquad &X^{\alpha_j}v_w=q^{-2}v_w, \cr
T_iv_w = qv_w, 
&\qquad &T_jv_w=-q^{-1}v_w.\cr
}
$$
For each of these cases one checks the braid relations
by  direct computation.  For type $G_2$ case (2) the calculations
can be simplified by observing that $T_1T_2T_1=T_2T_1T_2$ as
operators on $\CC$-span$\{v_w,v_{s_iw}\}$.  From this it follows that
$T_1T_2T_1T_2T_1T_2=T_2T_1T_2T_1T_2T_1$ as operators on
$\CC$-span$\{v_w,v_{s_iw}\}$.
\endpf

\medskip
\section 4. Calibrated Representations

\medskip

A finite dimensional $\tilde H$-module  
$$\hbox{$M$ is {\it calibrated} if $M_t^{\rm gen}=M_t$,\enspace 
for all $t\in T$.}$$  
A calibrated module $M$ is {\it really calibrated} if $t=t_r$ for all
$t\in T$, i.e. $t_c=1$ in the polar decomposition (2.1).   

Suppose that $t\in T$ is regular, i.e. the stabilizer $W_t$ of $t$ under the
action of $W$ is trivial.  Then $M(t) = \bigoplus_{w\in W} M_{wt},$
since each $M_{wt}^{\rm gen}$ is one dimensional and $M_{wt}\ne 0$ whenever
$M_{wt}^{\rm gen}\ne 0$.  Thus $M(t)$ is calibrated when $t$ is
regular.  So any quotient of $M(t)$ is calibrated and, by Proposition
2.6b, any irreducible $\tilde H$-module $M$ with regular central character is
calibrated.

\vfill\eject

\subsection{Classification of irreducible calibrated modules}

We shall eventually show that the modules $\tilde H^{(t,J)}$ constructed in
Theorem 3.1 are all the irreducible calibrated $\tilde H$-modules.  The following
Proposition shows that the formulas which define the
$\tilde H$-modules in Theorem 3.1 are more or less forced.

\prop  Let $M$ be a calibrated
$\tilde H$-module and assume that for all $t\in T$ such that $M_t\ne 0$,
$$\hbox{(A1)\quad $t(X^\alpha_i)\ne 1$ for all $1\le i\le n$,
\qquad and
\qquad (A2)\quad $\dim(M_t)=1$.}
$$ 
For each $b\in T$ such that $M_b\ne 0$ let $v_b$ be a nonzero vector in
$M_b$.  The vectors $\{v_b\}$ form a basis of $M$.
Let $(T_i)_{cb}\in \CC$ and $b(X^\lambda)\in \CC$ be given by
$$T_iv_b = \sum_c (T_i)_{cb}v_c  
\qquad\hbox{and}\qquad 
X^\lambda v_b = b(X^\lambda) v_b.$$
Then
\smallskip
\itemitem{(a)} $\displaystyle{
(T_i)_{bb} = {q-q^{-1}\over 1-b(X^{-\alpha_i}) }
},$
for all $v_b$ in the basis, 
\smallskip
\itemitem{(b)} If $(T_i)_{cb}\ne 0$ then $c=s_ib$,
\smallskip
\itemitem{(c)} $(T_i)_{b,s_ib}(T_i)_{s_ib,b} 
= (q^{-1}+(T_i)_{bb})(q^{-1}+(T_i)_{s_ib,s_ib}).$
\endprop
\pf
The defining equation for $\tilde H$,
$$
X^\lambda T_i - T_i X^{s_i\lambda} =
(q-q^{-1}) { X^\lambda-X^{s_i\lambda}\over 1-X^{-\alpha_i} },
$$
forces
$$
\sum_c \left( c(X^\lambda)(T_i)_{cb} - (T_i)_{cb}b(X^{s_i\lambda})
\right)v_c
=
(q-q^{-1}) { b(X^\lambda)-b(X^{s_i\lambda})\over 1-b(X^{-\alpha_i}) }v_b
$$
Comparing coefficients gives
$$\eqalign{
c(X^\lambda)(T_i)_{cb}-(T_i)_{cb} b(X^{s_i\lambda}) 
&= 0,
\qquad\hbox{if $b\ne c$, and } \cr
\cr
b(X^\lambda)(T_i)_{bb}-(T_i)_{bb}b(X^{s_i\lambda})
&=(q-q^{-1}){ b(X^\lambda)-b(X^{s_i\lambda})\over 1-b(X^{-\alpha_i}) }. \cr
}
$$
These relations give:
$$
\hbox{If \quad $(T_i)_{cb}\ne 0$ \quad then \quad $b(X^{s_i\lambda}) =
c(X^\lambda)$
 \quad for all $X^\lambda\in X$, and}
$$
$$
(T_i)_{bb} = {q-q^{-1}\over 1-b(X^{-\alpha_i}) }
\quad \hbox{ if $b(X^{-\alpha_i})\ne 1$ and $b(X^\lambda)\ne b(X^{s_i\lambda})$ for
some $X^\lambda\in X$.}
$$
By assumption (A1), $b(X^{\alpha_i})\ne 1$ for all $i$.
For each fundamental weight $\omega_i$, $X^{\omega_i}\in X$ and
$b(X^{s_i\omega_i})=b(X^{\omega_i-\alpha_i})\ne b(X^{\omega_i})$ since
$b(X^{\alpha_i})\ne 1$.  Thus we conclude that 
$$
T_i v_b = (T_i)_{bb} v_b + (T_i)_{s_ib,b}v_{s_ib},
\qquad\hbox{with}\qquad
(T_i)_{bb} = {q-q^{-1}\over 1-b(X^{-\alpha_i}) }.
$$
This completes the proof of (a) and (b).
By the definition of $\tilde H$ the vector
$$
T_i^2 v_b = ((T_i)_{bb}^2 +(T_i)_{b,s_ib}(T_i)_{s_ib,b})v_b
+((T_i)_{bb}+(T_i)_{s_ib,s_ib})(T_i)_{s_ib,b}v_{s_ib}
$$
must equal
$$
((q-q^{-1})T_i+1)v_b 
=((q-q^{-1})(T_i)_{bb} +1)v_b + (q-q^{-1})(T_i)_{s_ib,b}v_{s_ib}.
$$
Using the formula for $(T_i)_{bb}$ and $(T_i)_{s_ib,s_ib}$ we find
$(T_i)_{bb}+(T_i)_{s_ib,s_ib} = (q-q^{-1})$.  
So, by comparing coefficients of $v_b$, we obtain the equation
$$
(T_i)_{b,s_ib}(T_i)_{s_ib,b} = (q-(T_i)_{bb})((T_i)_{bb}+q^{-1})
= (q^{-1}+(T_i)_{bb})(q^{-1}+(T_i)_{s_ib,s_ib}).\qquad\hbox{\qed}
$$

\prop  Let $M$ be an irreducible
calibrated module. Then, for all $t\in T$ such that $M_t\ne 0$,
\smallskip
\item{(a)} $t(X^{\alpha_i})\ne 1$ for all $1\le i\le n$, and 
\smallskip
\item{(b)} $\dim(M_t)=1$.
\pf
(a)  The proof is by contradiction.
Assume that $t(X^{\alpha_i})=1$.  
Let $\tilde HA_1$ be the subalgebra of $\tilde H$ generated by
$T_i$ and $X^{\alpha_i}$ and view $M$ as an $\tilde HA_1$-module 
by restriction.  Let $m_t$ be a nonzero element of $M_t$.
There is an $\tilde HA_1$-module  homomorphism
$$
\matrix{
\phi \colon &M(t) &\longrightarrow &M \cr
&v_t &\longmapsto &m_t \cr
}
$$
where $M(t)$ is the (two dimensional) principal series module for 
$\tilde HA_1$ and $v_t$ is the generator of $M(t)$.
It is easy to check that when $t(X^{\alpha_i})=1$ the module $M(t)$
is an irreducible $\tilde HA_1$-module.
Thus the map $\phi$ is injective and we can view 
$M(t)$ as a submodule of $M$.  
A direct check shows that $M(t)$ is not calibrated and thus
it follows that $M$ is not calibrated. 
This is a contradiction to the assumption
that $M$ is calibrated.  Thus $t(X^{\alpha_i})\ne 1$.

(b)  
The proof is by contradiction.  Assume that $t\in T$ is such that
$\dim(M_t)>1$.  Let $m_t$ be a nonzero element of $M_t$.  Since $M$ is
calibrated, the action of any $\tau_i$ on any weight vector $m$ is a
linear combination of the action of $T_i$ and a multiple of the identity.
Thus, since $M$ is irreducible, we must be able to generate the rest of
$M_t$ by applying $\tau$-operators to $m_t$.  Since $\dim(M_t)>1$
there must be a sequence of $\tau$-operators such that 
$$n_t=\tau_{i_1}\tau_{i_2}\cdots \tau_{i_p}m_t$$
is a nonzero vector in $M_t$ which is not a multiple of $m_t$.
Assume that the sequence $\tau_{i_1}\tau_{i_2}\cdots \tau_{i_p}$
is chosen so that $p$ is minimal.

Let us defer, momentarily, the proof of the following claim.
\medskip\noindent
\hbox{{\it Claim:}  The element $s_{i_1}s_{i_2}\cdots s_{i_p}=1$ in $W$.}
\medskip\noindent
The claim implies that there is some
$1<k\le p$ such that
$s_{i_1}s_{i_2}\cdots s_{i_k}$ is not reduced and we can use the braid
relations to rewrite this word as
$s_{i_1'}\cdots s_{i_{k-2}'}s_{i_k}s_{i_k}$.
By Proposition 2.7e the $\tau_i$ operators also satisfy the braid relations and so 
$$n_t=\tau_{i'_1}\tau_{i'_2}\cdots 
\tau_{i'_{k-2}}\tau_{i_k}\tau_{i_k}\cdots \tau_{i_p}m_t.$$
By Proposition 2.7b, the operator $\tau_{i_k}\tau_{i_k}$ is equal to a constant
times the identity map and thus
$$n_t=c\tau_{i'_1}\tau_{i'_2}\cdots 
\tau_{i'_{k-2}}\tau_{i_{k+1}}\cdots \tau_{i_p}m_t,$$
where $c$ is some constant.  The constant $c$ must be nonzero
since $n_t$ is not $0$.  But the expression
$$c^{-1}n_t= \tau_{i'_1}\tau_{i'_2}\cdots 
\tau_{i'_{k-2}}\tau_{i_{k+1}}\cdots \tau_{i_p}m_t$$
is shorter than the original expression of $n_t$ and this contradicts 
the minimality of $p$.
It follows that $\dim(M_t)\le 1$.

\smallskip\noindent
{\it Proof of the claim.}
By [St, 3.15, 4.2, 5.3] the stabilizer $W_t$ of $t$ under the action
of $W$ is
$$W_t =  \langle s_\alpha \ | \ \alpha\in Z(t) \rangle
\qquad\hbox{where}\qquad Z(t)=\{\alpha>0\ |\ t(X^\alpha)=1\}.$$
The elements of the orbit $Wt$ can be identified with the
cosets of $W/W_t$ and these can be identified with the chambers of
$\RR^n\setminus (\bigcup_\alpha H_\alpha)$ which are on the positive
side of all the hyperplanes $H_\alpha$ for $\alpha\in Z(t)$.
Specifically, the element $t\in Wt$ corresponds to the chamber
$C$ and the element $wt$ of $Wt$ corresponds to the chamber $w^{-1}C$.

For any $1\le j\le p$ we have that
$(s_{i_{j+1}}\cdots s_{i_p}t)(X^{\alpha_{i_j}})\ne 1$,
since $\tau_{i_j}\cdots \tau_{i_p}m_t$ is well defined.  This means
that $s_{i_j}(s_{i_{j+1}}\cdots s_{i_p}t)(X^{\omega_{i_j}})
=(s_{i_{j+1}}\cdots s_{i_p}t)(X^{\omega_{i_j}-\alpha_{i_j}})\ne
(s_{i_{j+1}}\cdots$ $s_{i_p}t)(X^{\omega_{i_j}})$ and thus that
$s_{i_j}(s_{i_{j+1}}\cdots s_{i_p}t)\ne s_{i_{j+1}}\cdots s_{i_p}t$. 
So $s_{i_j}\cdots s_{i_p}t$ and $s_{i_{j+1}}\cdots s_{i_p}t$ both correspond to
chambers on the positive side of all the hyperplanes $H_\alpha$,
$\alpha\in Z(t)$.  These two chambers have a common face and this face
is contained in the hyperplane
$H_{s_{i_p}\ldots s_{i_{j+1}}\alpha_j}$.

In this way we can identify the sequence
$t, s_{i_p}t, s_{i_{p-1}}s_{i_p}t,\ldots, s_{i_1}\cdots s_{i_p}t$
with a sequence of chambers where successive chambers in the sequence
are adjacent (share a face) and all the chambers in the sequence
are on the positive side of all the hyperplanes $H_\alpha$, $\alpha\in Z(t)$.
Since $s_{i_1}\ldots s_{i_p}t=t$, the first and last chamber in this sequence
are the same.  It follows that $s_{i_1}\cdots s_{i_p}=1$ in $W$.
\endpf

\prop  Let $M$ be an irreducible calibrated $\tilde H$-module.
Suppose that $M_t$ and $M_{s_it}$ are both nonzero.  Then
the map $\tau_i\colon M_t\to M_{s_it}$ is a bijection.
\pf  By Proposition 4.2b, $\dim(M_t)=\dim(M_{s_it})=1$, and
thus it is sufficient to show that $\tau_i$ is not the zero map.  Let $v_t$ be a
nonzero vector in $M_t$.   Since $M$ is irreducible there must be a sequence of
$\tau$ operators such that 
$$v_{s_it}=\tau_{i_1}\cdots \tau_{i_p}v_t$$
is a nonzero element of $M_{s_it}$.  Let $p$ be minimal such that
this is the case.
We have $\tau_i\tau_{i_1}\cdots \tau_{i_p}v_t\in M_t$.
Using the claim which was proved in the proof of Proposition 4.2 we have that
$s_is_{i_1}\cdots s_{i_p}=1$ in $W$. For notational convenience $i_0=i$. Let
$0\le k< p$ be maximal such that $s_{i_k}s_{i_{k+1}}\cdots s_{i_p}$ is not
reduced.  If
$k\ne 0$ then we can use the braid relations to get
$$v_{s_it}=\tau_{i_1}\cdots\tau_{i_k}\tau_{i_k}
\tau_{i'_{k+2}}\cdots\tau_{i'_p}v_t.$$
By Proposition 2.7c $\tau_{i_k}\tau_{i_k}$ is a multiple of the identity
and so
$$v_{s_it}=c \tau_{i_1}\cdots\tau_{i_{k-1}}
\tau_{i'_{k+2}}\cdots\tau_{i'_p}v_t.$$
But this contradicts the minimality of $p$.  Thus we must have
$k=0$, $p=1$ and
$$v_{s_it}=\tau_i v_t.$$
Thus, since $v_{s_it}\ne 0$, $\tau_i\ne 0$.  
\endpf

\prop If $M$ is a calibrated $\tilde H$-module and $t\in T$ is
such that $M_t\ne 0$ then $t$ is calibratable for all
$R_{ij}$ generated by simple roots $\alpha_i$ and $\alpha_j$
in $R$.
\pf   Let $\tilde H_{ij}$
be the subalgebra of $\tilde H$ generated by $T_i, T_j, X^{\alpha_i},$
and $X^{\alpha_j}$ and view $M$ as an $H_{ij}$ module by restriction.
The irreducible representations of rank two affine
Hecke algebras have been classified and constructed explicitly
in [Ra3].  From this classification it is easy to check that
the only weights $t$ which appear in calibrated
$\tilde H_{ij}$-modules are those that are calibratable
for $R_{ij}$.  Thus, if $M_t\ne 0$, then $t$ must be calibratable
for $R_{ij}$. 
\endpf

\thm  Let $M$ be an irreducible calibrated $\tilde H$-module.
Let $t$ be the central character of $M$ and let 
$J=R(w)\cap P(t)$ for any $w\in W$ such that $M_{wt}\neq 0$.
Then $(t,J)$ is a placed skew shape and 
$$M\cong \tilde H^{(t,J)},$$ 
where $H^{(t,J)}$ is the module defined in Theorem 3.1.
\endthm
\pf
Proposition 4.3 shows that if $M_t$ and $M_{s_it}$ are both nonzero
then both $\tau_i\colon M_t\to M_{s_it}$ and $\tau_i\colon M_{s_it}\to M_t$
are bijections.  Thus, by Proposition 2.7d $t(X^{\alpha_i})\ne q^{\pm 2}$
and so there is an edge $t\leftrightarrow s_it$ in the calibration graph.
This shows that $\supp(M)$ is a single connected component of the calibration
graph $\Gamma(t)$.  Then $(t,J)$ (as defined in the statement of the
Theorem) is the corresponding placed shape. By Proposition 4.4 $(t,J)$ must be a
placed {\it skew} shape. Propositions 4.1 and 4.2 show that there is
at most one calibrated $\tilde H$-module $M$ such that $\supp(M)$ is the connected
component of $\Gamma(t)$ labeled by $(t,J)$. Thus we must have that
$M\cong \tilde H^{(t,J)}$.
\endpf

\vfill\eject

\centerline{\smallcaps References}

\medskip
\item{[A]} {\smallcaps S.\ Ariki}, 
{\it  On the decomposition numbers of the Hecke algebra of $G(m,1,n)$}, 
J.\  Math.\  Kyoto Univ.\ {\bf 36} (1996), 789--808.

\medskip
\item{[AK]} {\smallcaps S.\ Ariki and K.\ Koike}, {\it A Hecke algebra of
$(\ZZ/r\ZZ)\wr S_n$ and construction of its irreducible representations},
Adv.\ in Math.\ {\bf 106} (1994), 216--243.


\medskip
\item{[B]} {\smallcaps N.\ Bourbaki}, 
{\it Groupes et alg\`ebres de Lie, Chapitres 4,5 et 6}, 
Elements de Math\'e\-matique, Hermann, Paris 1968.


\medskip
\item{[CG]} {\smallcaps N.\ Chriss and V.\ Ginzburg}, 
{\sl Representation Theory and Complex Geometry}, Birkh\"auser, 1997.

\medskip
\item{[H]} {\smallcaps P.N.\ Hoefsmit}, 
{\it Representations of Hecke algebras of finite groups with $BN$-pairs 
of classical type}, 
Ph.D.\ Thesis, University of British Columbia, 1974.

\medskip
\item{[HO1]} {\smallcaps G.J.\ Heckman and E.M.\ Opdam}, 
{\it Yang's system of particles and Hecke algebras}, 
Ann. of Math. (2) {\bf 145} (1997), 139--173.

\medskip
\item{[HO2]} {\smallcaps G.J.\ Heckman and E.M.\ Opdam}, 
{\it Harmonic analysis for affine Hecke algebras}, Current Developments
in Mathematics, Intern. Press, Boston 1996.

\medskip
\item{[IM]} {\smallcaps N.\ Iwahori and H.\ Matsumoto}, 
{\it On some Bruhat decomposition and the structure of the Hecke rings of 
$p$-adic Chevalley groups}, 
Publ.\ Math.\ IHES {\bf 40} (1972), 81--116.

\medskip
\item{[Ka]} {\smallcaps S-i.\ Kato},
{\it Irreducibility of principal series representations for Hecke
algebras of affine type}, 
J. Fac. Sci. Univ. Tokyo Sect. IA Math. {\bf 28} (1981), 929--943.

\medskip
\item{[KL]} {\smallcaps D.\ Kazhdan and G.\ Lusztig}, 
{\it Proof of the Deligne-Langlands conjecture for Hecke algebras}, 
Invent. Math. {\bf 87} (1987), 153--215.

\medskip
\item{[Lu1]} {\smallcaps G.\ Lusztig},
{\it Singularities, character formulas, and a $q$-analog of weight
multiplicities}, Analysis and topology on singular spaces, II, III (Luminy, 1981),
 Ast\'erisque {\bf 101-102}, Soc. Math. France, Paris, 1983, 208--229.


\medskip
\item{[Lu2]} {\smallcaps G.\ Lusztig},
{\it Affine Hecke algebras and their graded version}, J. Amer. Math.
Soc. {\bf 2} (1989), 599--635.

\medskip
\item{[M]} {\smallcaps I.G.\ Macdonald}, 
{\it Affine Hecke algebras and orthogonal polynomials}, 
S\'eminaire Bourbaki, 47\`eme ann\'ee, ${\rm n}^{\rm o}$ 797, 1994--95,
Ast\'erisque {\bf 237} (1996), 189--207.

\medskip
\item{[Ma]} {\smallcaps H.\ Matsumoto}, 
{\sl Analyse harmonique dans les syst\`emes de Tits bornologiques de type affine},
Lect. Notes in Math. {\bf 590}, Springer-Verlag, Berlin-New York,
1977.


\medskip
\item{[Ra1]} {\smallcaps A.\ Ram}, 
{\it Seminormal representations of Weyl groups and Iwahori-Hecke algebras}, 
Proc.\ London Math.\ Soc.\ (3) {\bf 75} (1997), 99-133.

\medskip
\item{[Ra2]} {\smallcaps A. \ Ram}, 
{\it Standard Young tableaux for finite root systems}, preprint 1998.

\medskip
\item{[Ra3]} {\smallcaps A.\ Ram}, 
{\it Irreducible representations
of rank two affine Hecke algebras}, preprint 1998.

\medskip
\item{[Ra4]} {\smallcaps A.\ Ram}, 
{\it Skew shape representations are irreducible}, preprint 1998.

\medskip
\item{[RR1]} {\smallcaps A.\ Ram and J. Ramagge}, 
{\it  Jucys-Murphy elements come from affine Hecke algebras}, in preparation.

\medskip
\item{[RR2]} {\smallcaps A.\ Ram and J. Ramagge}, {\it Calibrated
representations and the $q$-Springer correspondence}, in preparation.

\medskip
\item{[Ro]} {\smallcaps F.\ Rodier}, 
{\it D\'ecomposition de la s\'erie principale des groupes r\'eductifs $p$-adiques},
Noncommutative harmonic analysis and Lie groups (Marseille, 1980),
Lect. Notes in Math. {\bf 880}, Springer, Berlin-New York 1981, 408--424. 

\medskip
\item{[St]} {\smallcaps R.\ Steinberg}, 
{\it Endomorphisms of linear algebraic groups}, Mem. Amer.
Math. Soc. {\bf 80}, Amer. Math. Soc., Providence, R.I. 1968. 

\medskip
\item{[Wz]}  {\smallcaps H. Wenzl}, 
{\it Hecke algebras of type $A\sb n$ and subfactors}, Invent. Math. {\bf 92}
(1988), 349--383.

\vfill\eject
\end

  The lattices
$$Q = \sum_{i=1}^n \ZZ\alpha_i 
\qquad\hbox{and}\qquad
P= \sum_{i=1}^n \ZZ\omega_i$$
are the {\it root lattice} and the {\it weight lattice} respectively.
Let $L$ be a $W$-invariant lattice such that $Q\subseteq L\subseteq P$.

Assume that $q\in \RR_{>0}$ and $q\ne 1$.
This is a minor assumption since the 
representation theory of $\tilde H$ will be ``the same'' 
for all $q\in \CC^*$ that are not roots of unity. 
An element $t\in T$ is {\it real} if $t(X^\alpha)\in \RR_{>0}$ for 
all $\alpha\in R$.  If $t$ is real then
$$\matrix{
t\colon &X(T) &\longrightarrow &\CC^* \cr
&X^\lambda &\longmapsto &q^{2\langle \gamma,\lambda\rangle} \cr
}$$
for an appropriate $\gamma\in \RR^n$.  
A {\it real central character} is the $W$-orbit of a real $t\in T$.

The advantage of working with real $t\in T$ is that there is always
a convenient representative of the $W$-orbit of $t$.  An element
$\gamma\in \RR^n$ is {\it dominant} if $\langle \gamma,\alpha\rangle\ge 0$
for all $\alpha\in R$.  The dominant $\gamma\in \RR^n$ index the real 
central characters.

If $\gamma\in \RR^n$ define
$$Z(\gamma) = \{ \alpha\in R^+\ |\ \langle\gamma,\alpha\rangle=0 \}
\qquad\hbox{and}\qquad
P(\gamma)= \{ \alpha\in R \ |\ \langle\gamma,\alpha\rangle=1 \}.
$$
A {\it placed shape} is a pair $(\gamma,J)$ such that $\gamma\in \RR^n$ is
dominant and $J\subseteq P(\gamma)$.  If $(\gamma,J)$ is a placed shape define
$$\cF^{(\gamma,J)}=\{ w\in W \ | \ R(w)\cap Z(\gamma)=\emptyset,
\ R(w)\cap P(\gamma) = J \},$$
where, for an element $w\in W$, the set
$R(w)=\{ \alpha>0 \ |\ w\alpha < 0\}$
is the {\it inversion set} of $w$.
The elements of the set $\cF^{(\gamma,J)}$ are called {\it standard Young
tableaux} of shape $(\gamma,J)$.
The following is one of the main results of [Ra1].

\thm  Let $\gamma\in \RR^n$ be dominant and let $t$ be the element of $T$
given by $t(X^\lambda) = q^{2\langle \gamma,\lambda\rangle}$ for all
$X^\lambda\in X(T)$.  Then
the connected components of the calibration graph $\Gamma(t)$ are given by 
the partition of the vertices according to the sets
$$\cF^{(\gamma,J)}t,
\qquad\hbox{such that $J\subset P(\gamma)$ and $\cF^{(\gamma,J)}\ne \emptyset$.}
$$
\endthm